\documentclass[conference,10pt]{IEEEtran}
\IEEEoverridecommandlockouts
\usepackage{cite}
\usepackage{amsmath,amssymb,amsfonts}
\usepackage{graphicx}
\usepackage{textcomp}
\usepackage{xcolor}
\def\BibTeX{{\rm B\kern-.05em{\sc i\kern-.025em b}\kern-.08em
    T\kern-.1667em\lower.7ex\hbox{E}\kern-.125emX}}

\usepackage{algorithm}
\usepackage[noend]{algpseudocode}

\newtheorem{example}{Example}

\usepackage{tikz}

\usetikzlibrary{arrows,automata,backgrounds,fit,shapes,decorations.pathreplacing}

\RequirePackage[T1]{fontenc}
\RequirePackage[utf8]{inputenc}

\usepackage{amsmath,amssymb}
\usepackage{graphicx}

\newcommand{\boxtensor}{{\Box\kern-9.03pt\raise1.42pt\hbox{$\times$}}}

\newcommand{\be}{\begin{eqnarray}}
\newcommand{\ee}{\end{eqnarray}}

\usepackage{fancyhdr}

\makeatletter
\newcommand*\titleheader[1]{\gdef\@titleheader{#1}}
\AtBeginDocument{%
  \let\st@red@title\@title
  \def\@title{%
    \bgroup\normalfont\large\raggedleft\@titleheader\par\egroup
    \vskip0.5em\st@red@title}
}
\makeatother

\begin{document}

%\IEEEoverridecommandlockouts
%\IEEEpubid{\makebox[\columnwidth]{978-1-7281-7116-6/20/\$31.00 ©2020 IEEE \hfill} \hspace{\columnsep}\makebox[\columnwidth]{ }}

\title{Comparison of Randomized Solutions for Constrained Vehicle Routing Problem}
%\titleheader{\small \textit{Proc. of the $2nd$ International Conference on Electrical, Communication and Computer Engineering (ICECCE)\\ 12-13 June 2020, Istanbul, Turkey}}
%\fancyhdr[C]{2016 IEEE 24th International Requirements Engineering Conference} %zentrierte Kopfzeile
\author{\IEEEauthorblockN{İbrahim Ethem Demirci}
\IEEEauthorblockA{Department of Mathematics,\\ Hacettepe University\\ Ankara, Turkey\\
demirciibrahim.e@gmail.com}
\and
\IEEEauthorblockN{Şaziye Ece Özdemir}
\IEEEauthorblockA{Department of Mathematics,\\Pamukkale University,\\ Denizli, Turkey and\\Department of Mathematics,\\  Hacettepe University,\\ Ankara, Turkey\\
saziyeo@pau.edu.tr}
\and
\IEEEauthorblockN{Oğuz Yayla}
\IEEEauthorblockA{Department of Mathematics,\\ Hacettepe University,\\Ankara, Turkey\\
oguz.yayla@hacettepe.edu.tr}
}

\maketitle
\IEEEpubidadjcol

\begin{abstract}
In this short paper, we study capacity constrained vehicle routing problem (CVRP) and its solution by randomized Monte Carlo methods. For solving CVRP we use some pseudorandom number generators commonly used in practice. We use linear, multiple-recursive, inversive and explicit inversive congruential generators and obtain random numbers from each to provide a route for CVRP. Then we compare the performance of pseudorandom number generators with respect to the total time the random route takes. We also constructed an open-source library github.com/iedmrc/binary-cws-mcs on solving CVRP by Monte-Carlo based heuristic methods.
%\subclass{MSC 05B10 \and MSC 94A55}
\end{abstract}


\section{Introduction}
The Monte Carlo method is used in computational algorithms where repeated random sampling is needed. It provides randomness to solve problems that might be deterministic in principle. Many of the most useful techniques use deterministic pseudorandom sequences because it is simple to test and re-run the simulation. Monte Carlo method has a lot of areas of application, for example, physical science, engineering, computational biology, applied statistics, computer graphics, finance, business and mathematics \cite{ref_proc1}.

Monte Carlo method  uses random samples; hence needs a good random number generation. There are two main types of  random number generation: true random generation and pseudorandom generation. We use a pseudorandom generator because of its deterministic behaviour, that is, one can re-run the generation process.

There are various pseudorandom generators in literature. In this study, we examine some of them: linear congruential method (LCG), inversive congruential generator (ICG), Lehmer congruential generation (LehmerCG) and  explicit inversive congruential generator (EICG). Then we test their performance on the capacitated vehicle routing problem (CVRP) on TSPLIB data sets\cite{TSPLIB}. %After all is said and done; aim of this paper obtain randomness results. //

\section{Related Work}
The vehicle routing problem (VRP) was first studied in \cite{DR1959} by George Dantzig and John Ramser in 1959, where petrol deliveries were considered algorithmically. The ultimate aim of the VRP is to minimize time, distance or budget cost of route. For other variants of VRP, see the recent study \cite{VLM2019}.

Several methods have been used for solving CVRP. There are heuristic and meta-heuristic methods. Heuristic methods speed up duration for reaching a reasonable and optimal solution, and they are based on four principles: understand the problem, make a plan, apply the plan, evaluate and adjust respectively. For example; Clark \& Wright's Savings \cite{CW1964} is used in heuristic methods thoroughly. See also its generalizations in \cite{AO2005,DC2011,S2006} and recent variants in
\cite{BSK2020,AGS2019}.

In 1964, Clarke and Wright developed a saving algorithm based on Dantzig and Ramser's greedy approximation, which is one of the most known heuristic solution methods for VRP. This is called as CWS algorithm. %Vehicles numbers are not constant this method, n vehicles serves one customer. Their algorithm follows a path like this: Firstly all pairs of nodes i,j  calculate the savings merging delivery locations and separate the savings in decreasing order. Starting at the top of the remaining list of saving, combination the two routes associated with the largest saving provided that: two delivery locations are not already on the same routes, the maximum capacity of vehicle is not surpassed, i and j are first or last customer on cycles. Repeat last step until no additional savings can be achieved.
In this study we consider Capacitated VRP (CVRP), which  is like VRP with the additional constrain that every vehicle has a capacity of the amount that it can carry. %On the other hand; the vehicle routing problem with time windows deals with the time windows every vehicle can be run in period of time and/or every customer has to be supplied among a certain period of time.

There are some studies on applying the Monte Carlo method to the VRP and CVRP.
For VRP it was applied the first time in the 1970s and random sampling for the distance constrained VRP applied in \cite{ref_lncs1}. %Monte Carlo Techniques are often used classic heuristic methods collapse because it is difficult to derive an acceptable assessment function to determine where the nominee successors is  selected later. These techniques are currently widely used in all kinds of optimization algorithms and  also be applied to the Vehicle Routing Problem.
Authors in \cite{ref_article1} use the Monte Carlo method in the CWS algorithm to provide a better solution for VRP. They also applied the Monte Carlo method for CVRP. In this paper we mainly based on the methods given in \cite{ref_article1} and their improvements. There are other studies combine Monte-Carlo simulation and heuristic method for CVRP, see \cite{Juan2011MC,Juan2014Simheuristic,Juan2019Random,Juan2019sim}.

On the other hand, metaheuristic methods has made great progress in the past 15 years for VRP  and classified as local search, population search, and learning algorithms, which were used much earlier than CVRP \cite{ref_article2}.

Another randomized methodology applied in solving VRP is GRASP (greedy randomized adaptive search procedure). It is a recursive metaheuristic and each recursive compose two stages that are construction and local search. It was for the first time exploited by Feo and Resende in 1989 and largely used in multiple applications \cite{ref_book1}. GRASP solutions might be generated in the same number of time need for the local search procedure to converge from a single random start.
%However, none of these studies mentioned about which randomness methods are using for choosing candidate element.

\section{Methodology}
The goal of this section is to describe the methodology we follow to provide a randomized solution for CVRP. Initially, there are a lot of pseudorandom numbers generation methods, see \cite{L1988} and \cite{L1999}. In this paper we consider only inversive congruential generator (ICG), multiple-recursive  generator (MRG), Lehmer congruential generator (LehmerCG), linear congruential generator (LCG). The main aim of this paper is to compare the CVRP solutions obtained by using these pseudorandom generation methods.

The first pseudorandom number generation method is Lehmer congruential generator (LehmerCG), also known as Park Miller random number generator, the general formula is $x_{k+1}=a \cdot x_{k} \mod m$ where  $m$ is a prime number or power of a prime number, $a$ is the multiplier and the seed $x_0$ is expected to be coprime to $m$.

The linear congruential generator (LCG), introduced in 1958 by W.E Thomson and A Rotenberg, is an extension of the Lehmer generator and it is the most used and studied generator. This method is still practically used since it is so simple to define. It has only two parameters namely a large integer $m$ and integers $a,b$ with gcd($a,m)=1$. Let $\mathbb{Z}_{m}=\{0,1,2,\ldots\,m-1\}$ and choose an initial value $x_{0}\in \mathbb{Z}_{m}$ such that gcd($x_{0},m)=1$. Then we generate a sequence $x_{n+1}=a \cdot x_{n} + b \mod m$ for $n=0,1, \ldots$ where $m$ is referred to as the modulus and $a$ is the multiplier.

%Taking into account all of these, we can create the algorithm as follows.
%\newline In the beginning; We choose a initial number and generate to  random numbers via linear congruential method sufficiently. Moreover this step maintain for the other seeds.

The third random number generation is multi-recursive generator (MRG). It generates a random sequence $x_{n}$  as $x_{n+1}=a \cdot x_{n} + b \cdot x_{n-1} \mod m$ for $n=1, \ldots$ and initial values $x_0,x_1 \in \mathbb{Z}_m$ where $m$ is referred to as the modulus and $a,b$ are multipliers.

The fourth random number generation is the inversive congruential generator (ICG), which is a nonlinear generator introduced by Eichenauer and Lehn in 1986. In this method, $p$ is a prime number, $x^{-1} \cdot x=1 \mod p$ where $x^{-1}$ is called the inverse of an integer $x$ and we generate random numbers with the standard formula for a chosen seed $x_{0}$ 
$$x_{i+1}=(a\cdot x_{i}^{-1}+c) \mod p \quad \mbox{if}\quad x_{i}\neq 0 $$
$$x_{i+1}=c \quad \mbox{if} \quad x_{i}=0. $$
%Therefore we generate random numbers and apply this method for other seeds.

Lastly, the explicit inversive congruential generator (EICG) was introduced by Eichenauer-Herrmann in 1993 and it is a different form of inversive generator. It is defined as $x_{k} = \overline{k + k_{0}} \mod p$ where $k_0$ is a seed value.

We use pseudorandom numbers at various algorithms such as Nearest Neighbor Insertion (NNI), Clarke $\&$ Wright's Savings Algorithm (CWS), Monte Carlo Simulation.

Below, we consider that we have $m$ vehicles with equal capacity $Q$ and they are assigned to $n$ jobs with demands $d_i$, $i=1,2,\ldots,n$.

\subsection{Nearest Neighbor Insertion}
\iffalse
Todo:\\
1. Lütfen ingilizce gramer ve kelime kullanım kurallarını gözden geçirip gerekli düzeltmeleri yapar misin? İngilizcemizi geliştirecegiz.\\
2. NNI ve CWS algoritmasini ekler misin? Algorithm lateks paketini kullanabilirsin. https://tr.overleaf.com/learn/latex/algorithms\\
3. Örnegi sonuna kadar yapalim ve butun araçlara isleri dagitalim.\\
4. NNI MMC ve algritması\\
5. CWS MMC ve algoritması\\
6. \\
\fi

The Nearest Neighbor Insertion (NNI) Algorithm is one of the methods for the vehicle routing problem. % which put forward useful application for solving problem. 
The NNI algorithm firstly selects a vehicle and a customer randomly  then adds another customer which is the nearest unvisited until the vehicle's capacity gets full. After that, we choose another vehicle, and similar steps are applied  until all customers are visited, or all vehicles' capacity is full, which depends on tightness. Tightness is a relation between all customers' demands and the total capacity of the vehicle, and defined as
\begin{center}
 $T=\dfrac{\sum d_i}{Qm}.$
\end{center}
If the T is very close to 1, the algorithm most likely does not yield a proper result.
We present the pseudocode  of the NNI method in  Algorithm 1.
\begin{example}
Let us examine the NNI method with an example where we have three vehicle $(m=3)$ with capacities $Q=8$, there are 8 customers $j_i$ $(i=1\ldots 8)$ with demands  $d_1=2, d_2=4, d_3=3, d_4=1, d_5=5, d_6=2, d_7=2, d_8=3$ respectively. In Table I, we give the distance between each customer.
\begin{table}[!h]
\center
\caption{Distance Matrix}
\begin{tabular} {|l|l|l|l|l|l|l|l|l|}
\hline
$j_i$ &$j_1$ &$j_2$ & $j_3$ & $j_4$ & $j_5$ & $j_6$ & $j_7$ & $j_8$\\
\hline
$j_1$ &0&10&20&13&12&7&15&9\\
\hline
$j_2$ &10&0&4&7&6&20&13&17\\
\hline
$j_3$ &20&4&0&11&6&19&14&8\\
\hline
$j_4$ &13&7&11&0&3&5&6&18\\
\hline
$j_5$ &12&6&6&3&0&9&12&16\\
\hline
$j_6$ &7&20&19&5&9&0&15&11\\
\hline
$j_7$ &15&13&14&6&12&15&0&21\\
\hline
$j_8$ &9&17&8&18&16&11&21&0\\
\hline
\end{tabular}
\end{table}
Firstly a customer $(j_1)$ is selected randomly, and NNI method chooses the nearest unvisited customer $(j_6)$, then checks whether the sum of their demands exceeds the capacity of the vehicle $(2+2\leq8)$. If enough capacity exists then the method constitutes a route as $[ j_1\quad j_6 ]$. If this is continued by selecting the nearest unvisited customer to $j_6$, route $[ j_1\quad j_6\quad j_4 ]$ is obtained for the first vehicle, if we add to $j_5$ for the nearest customer to $j_4$, this will surpass the capacity. Therefore pick $j_5$ for the second vehicle randomly and the one closest to it $j_2$ or $j_3$ but we prefer to $j_3$ because $j_2$ surpasses the capacity of the second vehicle thus the route is occurred as $[j_5\quad j_3]$. Lastly take the third vehicle and  pick one of the remaining jobs randomly $j_2$ and one of the nearest to it $j_7$ in the rest of jobs. Finally, the latest customer $j_8$ and its demand is suitable for the vehicle's capacity, so the route constituted for the third vehicle like $[j_2\quad j_7 \quad j_8]$.

\end{example}

\begin{algorithm}
\caption{Nearest Neighbor Insertion}
\begin{algorithmic}[1]
\Require{j, m, distance, demand, Q}
\State{capacity = Initialize to a zero array of length $m$}
\State{routes = Initialize to an empty two dimensional array}
\State v = 0
\While{v $<$ m}  
    \State $j_0$ = random unvisited job
    \While{capacity[v] + demand[$j_0$] $\leq$ Q}
    \State Append $j_0$ to routes[v]
    \State Update capacity[v] = capacity[v] + demand[$j_0$]
    \State $j_0$ = Nearest Unvisited Customer to $j_0$
    \EndWhile  \label{roy's loop}
    \State v = v + 1
\EndWhile
\State \Return{routes}
\end{algorithmic}
\end{algorithm}

\subsection{Clark \& Wright's Savings (CWS) Algorithm}
Clark and Wright \cite{CW1964} have published an algorithm for the vehicle routing problem based upon saving opinion, which is a heuristic algorithm, so the solution may not be provided with an optimal solution. In CWS algorithm, we have an ($n+1) \times (n+1)$ distance matrix $C$ whose entries $c_{ij}$ are the distance between two customers $i$ and $j$ for $i,j=0,1,2,\ldots,n$, where we set the index 0 to the starting point of the vehicles, say depot. Then, the saving list $s_{ij}$ is created by using the following formula
$$s_{ij}=c_{0i}+c_{0j}-c_{ij}$$ for $i,j=1,2,\ldots ,n$. The steps of the CWS method are given in Algorithm 2 where the second part between steps 6 and 12 is called the \textit{process} procedure of (i,j) pair, and this process part is used in the next section. Let us examine this algorithm with an example.
\begin{example}
 We have three vehicles (m=3) with capacity 10 (Q=10) and there are  eight jobs (j=8) with demands given in Table \ref{tab:demand}.
\begin{table}
\center
\caption{Customer Demand}
\label{tab:demand}
\begin{tabular} {|l|l|}
\hline
customer&demand\\
\hline
$j_1$&2\\
\hline
$j_2$&4\\
\hline
$j_3$&3\\
\hline
$j_4$&1\\
\hline
$j_5$&5\\
\hline
$j_6$&2\\
\hline
$j_7$&2\\
\hline
$j_8$&1\\
\hline
\end{tabular}
\end{table}
\begin{table}
\center
\caption{Distance Matrix}
\begin{tabular}{|l|l|l|l|l|l|l|l|l|l|}
\hline
$c_{ij}$ & $j_0$ &$j_1$&$j_2$&$j_3$&$j_4$&$j_5$&$j_6$&$j_7$&$j_8$\\
\hline
$j_0$ &0&10&20&13&12&7&15&9&10\\
\hline
$j_1$&10&0&4&7&6&20&13&17&9\\
\hline
$j_2$ &20&4&0&11&6&19&14&8&13\\
\hline
$j_3$ &13&7&11&0&3&5&6&18&14\\
\hline
$j_4$ &12&6&6&3&0&9&12&16&8\\
\hline
$j_5$&7&10&19&5&9&0&15&11&9\\
\hline
$j_6$ &15&13&14&6&12&15&0&21&13\\
\hline
$j_7$ &9&17&9&18&16&11&21&0&13\\
\hline
$j_8$&10&9&13&14&8&9&11&13&0\\
\hline
\end{tabular}
\end{table}
\begin{table}
\center
\caption{Symmetric Saving Matrix}
\begin{tabular} {|l|l|l|l|l|l|l|l|l|}
\hline
$s_{ij}$ & $j_1$&$j_2$&$j_3$&$j_4$&$j_5$&$j_6$&$j_7$&$j_8$\\
\hline
$j_1$&0&26&23&16&7&8&2&13\\
\hline
$j_2$& &0&22&26&8&21&20&17\\
\hline
$j_3$ & & &0&22&15&22&4&9\\
\hline
$j_4$ & & & &0&10&15&5&14\\
\hline
$j_5$ & & & & &0&7&5&8\\
\hline
$j_6$ & & & & & & 0&3&14\\
\hline
$j_7$ & & & & & & &0&6\\
\hline
$j_8$ & & & & & & & & 0\\
\hline
\end{tabular}
\end{table}
If we take into account saving matrix, constitutes saving list from the high savings pair to the lowest saving pair such that 1-2, 2-4, 1-3, 2-3, 3-4, 3-6, 2-6, 2-7, 2-8, 1-4, 3-5, 4-6, 4-8, 6-8, 4-5, 1-6, 3-8, 2-5, 1-5, 5-6, 7-8, 4-7, 6-7, 1-7 respectively. First, we choose to the highest saving pairs 1-2 and 2-4, then check whether the sum of their demands to surpass the capacity of a vehicle or not, $(2+4+1 \leq 10)$, if not, 1-2-4 will be added to one route, then 1-3 is added to a route because it does not exceed total capacity of the first vehicle. Therefore a route is created for the first vehicle such that 3-1-2-4. We can skip the pairs which are selected previously because no need to include them in another route. For the second vehicle, we first pick the highest saving 6-8. Next, we choose the pair 5-6 and add it to the route as 5-6-8 as total demand ($2+1+5\leq10$) ensures the second vehicle's capacity. The next pair is 7-8 but we cannot add to the route because it exceeds the total capacity. Thus the last vehicle will go to the last customer 7 only.
\end{example}

\begin{algorithm}
\caption{Clark \& Wright's Savings}
\begin{algorithmic}[1]
    \Require{distancematrix and the number of jobs $n$} 
    %\State s = empty saving list
    \For {i in range (1,n+1)}
    \For {j in range (1,i)}
    \State {$s_{ij}$=distancematrix(0,i)+distancematrix(0,j)-distancematrix(i,j)}
    \State Add [$i,j,s_{ij}$]  into $s$
    \EndFor
    \EndFor
    \State Sort $s$ descending with respect to $s_{ij}$
    \For{[$i,j,s_{ij}$] in $s$}
        \If {Neither i nor j assigned to  a route in  routelist}
        \State Initiative  a new route with  ($i,j$)
        \EndIf
        \If  {$i$ or $j$ exists at the end of a route}
        \State Link ($i,j$) is added  to  this route
        \EndIf
        \If {Both $i$ and $j$ exist at the end of the routes }
        \State Two routes merge into one route
        \EndIf
    \EndFor
\State \Return{routelist}
\end{algorithmic}
\end{algorithm}

\subsection{Monte Carlo Simulation}
The Monte Carlo Simulation is used in the heuristic methods. Firstly one chooses a root node  from the search space. We apply $r$ random simulation to its children nodes, by going from the current node until a leaf node. This algorithm visits $nr$ nodes at end of the each step, and so  has the highest score or highest average out of $r$ probes opt for as a "best child" and then the best child is used for expansion.

\subsection{MCS-NNI Method}
This method is used with a probability function as described in the following equation, instead of selecting the nearest unvisited customer.
\begin{center}
 $p_{v_i}=1-\dfrac{f(v_i)}{\bar{n}}$
\end{center}
where $f(v_i) \in [\bar{n},1]$ indicates the rank of customers such that the nearest customer has rank 1, another nearest customer's rank 2, etc. The step is used for making selections (it can be thought of as selecting n children) and then $r$ simulations applied to each of n children as a Monte Carlo Simulation. On the other hand, there is not much information about the NNI method because this method is outperformed by the CWS method in the literature.

\subsection{Binary-CWS and Binary-CWS-MCS Method}
This method uses Monte Carlo Simulation and Clarke $\&$ Wright's savings list. For each method, the saving list is ordered from top to bottom. Then generated random numbers are used for the election of the  current saving pair.
In the Binary-CWS  method, whose pseudocode is presented in Algorithm 3, we pick a $p$ probability value randomly in $0.05\leq p \leq 0.40$ as mentioned in \cite{ref_article1}. The most crucial point is that if $p$ value is very close to 0, the algorithm gives solution almost equal to CWS algorithm or larger value of $p$ will end up with chaos. Thus, we choose $p$ value such that $0.05\leq p \leq 0.40$.

\begin{algorithm}
\caption{Binary-CWS Algorithm}
\begin{algorithmic}[1]
    \Require{savinglist, probability}
    \State routelist = []
    \State newlist = [1,2,...,Len(savinglist)]
    \While {newlist is not empty}
        \For {i in newlist}
        \State rand = random() mod 100
        \If {rand $ \geq $ probability*100}
            \State process (savinglist(i), routelist)
            \State discard i from newlist
        \EndIf
        \EndFor
    \EndWhile
\State \Return score(routelist), route
\end{algorithmic}
\end{algorithm}
In the Binary-CWS method, we need saving-list from top to bottom and a probability value $p$  between $0.05\leq p \leq 0.40$. We first pick a pair for the first vehicle, which has the highest score in the saving-list. We note that we use the parallel version method in the CWS algorithm and create routes for many vehicles at the same time. We generate random numbers.  Then we compare $p$ value with the chosen random number  (mod 100),  whether to continue with the next pair or not. If random number (mod 100) $\geq 100p$, then we choose the next highest-scoring pair. When the demand is smaller than capacity, the pair is added to the route. Otherwise this job is  assigned to another suitable job. For the next steps of the algorithm we scan the saving-list, because we do not choose pair even it contains one of the common job  previously selected pairs. By repeating all these steps, a route is created for each vehicle. Let us explain with an example. 

\begin{example}
There are 4 vehicles (m=4), 12 jobs (n=12) also customers' demands are ($d_1$=1200, $d_2$=1700, $d_3$=1500, $d_4$=1400, $d_5$=1700, $d_6$=1400, $d_7$=1200, $d_8$=1900, $d_9$=1800, $d_{10}$=1600, $d_{11}$=1700, $d_{12}$=1100) respectively and distance matrix  given in Table V. Assume that vehicles have equal capacity (Q=6000). We compute the saving list 7-5, 10-5, 8-5, 10-7, 7-4, 11-5, 5-3, 5-4, 5-2, 7-2, 7-3, 12-7, 11-10, 12-5, 11-7, 8-7, 11-8,9-7, 11-4, 10-6, 7-6, 10-8, 12-10, 10-9, 11-3, 11-9, 10-4, 9-5, 12-9, 12-4, 12-11, 12-3, 11-2, 9-8, 8-6, $\ldots$. On the other hand, we generate random numbers with Lehmer Random Number Generator (seed: a=16807, z=172361, m=$2^{31}$-1) and these numbers (mod 100) are 61, 80, 52, 61, 47, 14, 12, 8, 35, 55, 61, 40, 96, 90, 36, 49, 30, 34, 47, 45, 92, 57, 25, 97, 87, 98, 26, 42, 88, 82, 56, 76, 55, 40, 80, 73, 35, 45, 64, 67, 91, 74, 1, 47, 93, 73, 29, 65, 68, 46, 37, 50, 59, 87, 83, 67, 37, 32, 16, 66, 97, 20, 6, 53, 28, 35, 59, 13, 2, 29, 96, 84, 24, 9, 11, 36, 10, 7, 92, 40, 18, 14, 58, 71, 47, 40, 19, 93, 40, 62, 93, 97, 44, 58, 6, 24, 68, 54, 36, 25,  $\ldots$ respectively. Also we take  $p$ value  randomly such as 0.30.
Consider the parallel version and select pair 7-5  for the first vehicle. The first  random number $61\geq 30$ so we select the pair. We generate random number $80\geq 30$ for the second pair 10-5  and check the total demand $d_{7}+d_{5}+d_{10}=4500$ less than capacity so we continue to select the pair but pass the  pair 8-5 as it surpasses the capacity of the vehicle.  The 4th job pair 10-7 is selected. But it is as connected with jobs 7 and 10. Finally, we choose the pair 7-4 , generate number $54\geq 30$ and $d_{7}+d_{5}+d_{10}+ d_{4}=5900 \leq 6000$. Hence a new route is created  $v_1$ = [4  7  5  10].
Then we scan the saving-list and pass to pair which included the same job with first vehicle routes 4-7-5-10. For example; 11-5 is connected with 7-5 (due to job 5) and 9-7 (due to job 7) so we skip it.  We choose the disconnected pair 11-8 and generated random number 61 is bigger than 30, $d_{11}+d_{8}=3600\leq 6000$ so 11-8 is selected for the second vehicle. We continue  scanning the saving list till 12-9, generating random number $47\geq 30$ but total demand $d_{11}+d_{8}+d_{12}+d_{9}=6500\geq 6000$ surpasses the capacity so we pass to the disconnected next pair 12-11. We don't select it because generating number 14 is less than 30. When similarly continued, pair 9-8 is selected, random number $35\geq 30$ moreover total demand is $d_{11}+d_{8}+d_{9}=5400$ smaller than the capacity of the second vehicle. Thereby route is completed for the second vehicle as $v_2$ = [11  8  9].
If all steps applied similarly, we obtain routes for vehicles $v_3$ = [3  2  12 6], $v_4$ = [1]  respectively. So we finish the example.
\end{example}
 
We now continue with the Binary-CWS-MCS method, which uses the Monte Carlo simulation. Let us explain the method by taking into consideration Example 3. In addition to all algorithmic steps, we apply the simulation (1000 Binary-CWS simulation) for both state of choosing and state of not choosing. At the end of these simulations, we calculate the average score and choose which average score is higher than the other. If we continue to apply each pair in this way, we will create a binary-tree. We give the steps of Binary-CWS-MCS method given in Algorithm 4.

\begin{table}
\tiny
\center
\caption{Distance Matrix}
\begin{tabular} {|l|l|l|l|l|l|l|l|l|l|l|l|l|}
\hline
  0 & 9 & 14 & 23 & 32 & 50 & 21 & 49 & 30 & 27 & 35 & 28 & 18\\
\hline
  9 & 0 & 21 & 22 & 36 & 52 & 24 & 51 & 36 & 37 & 41 & 30 & 20\\
\hline
  14 & 21 & 0 & 25 & 38 & 5 & 31 & 7 & 36 & 43 & 29 & 7 & 6\\
\hline
 23 & 22 & 25 & 0 & 42 & 12 & 35 & 17 & 44 & 31 & 31 & 11 & 6\\
\hline
 32 & 36 & 38 & 42 & 0 & 22 & 37 & 16 & 46 & 37 & 29 & 13 & 14\\
\hline
 50 & 52 & 5 & 12 & 22 & 0 & 41 & 23 & 10 & 39 & 9 & 17 & 16\\
\hline
 21 & 24 & 31 & 35 & 37 & 41 & 0 & 26 & 21 & 19 & 10 & 25 & 12\\
\hline
 49 & 51 & 7 & 17 & 16 & 23 & 26 & 0 & 30 & 28 & 16 & 27 & 12\\
\hline
 30 & 36 & 36 & 44 & 46 & 10 & 21 & 30 & 0 & 25 & 22 & 10 & 20\\
\hline
 27 & 37 & 43 & 31 & 37 & 39 & 19 & 28 & 25 & 0 & 20 & 16 & 8\\
\hline
 35 & 41 & 29 & 31 & 29 & 9 & 10 & 16 & 22 & 20 & 0 & 10 & 10\\
\hline
 28 & 30 & 7 & 11 & 13 & 17 & 25 &27 & 10 & 16 & 10 & 0 & 10\\
\hline
 18 & 20 & 6 & 6 & 14 & 16 & 12 & 12 & 20 & 8 & 10 & 10 & 0\\
\hline
\end{tabular}
\end{table}

\begin{algorithm}
\caption{Binary-CWS-MCS Algorithm}
\begin{algorithmic}[1]
    \Require{distancematrix}
    \State get saving-list (distance matrix,n)
    \State list-ordered = order list descending with respect to savings $s_{ij}$
    \State route list [       ]  
    \While{list-ordered is not empty for s in range}
        \State decide whether to "process" or "skip" s by Monte Carlo
        \State $t_1$ average of 1000 calls of score Binary-CWS(s,list-ordered)
        \State $t_2$ average of 1000 calls of score Binary-CWS(s+1,list-ordered)
        \If {$t_1\geq t_2$}
        \State process (list-ordered[s],route-list)
        \State discard list-ordered[s] from list ordered
        \State return score (route-list),route list
        \EndIf
    \EndWhile
\end{algorithmic}
\end{algorithm}

\section{Results}
In this section we compare the results of different pseudorandom generators with different seeds under the CVRP problem. Firstly we take different CVRP problems form TSPLib \cite{TSPLIB} such as E-13-k4, E-n22-k4, E-n23-k3, E-n30-k3, E-n31-k7, E-n33-k4, E-n51-k5, E-n76-k7, E-n76-k8, E76-k10. For instance E-n51-k5 means that there are 51 jobs and 5 vehicles.  
Then we calculate the costs of each route by using the pseudorandom numbers generators given in Table VI.  We get the shortest route-list as a result for each case. We present the scores in  Table VII.

\begin{table}[!t]
\small
\caption{Parameters of pseudorandom number generators (PRNG)}
\begin{tabular}{|l|l|l|l|}
\hline
PRNG &Seed\\

\hline
LehmerCG &$a:48271, z:172361, m:2^{31}-1$ \\

\hline
LCG &$ [a:6364136223846793005, z:[172361], m:2^{64}$ \\

\hline
MRG &$ [a:1071064 , z:[135623,172361], m:2^{31}-19$ \\

\hline
ICG &$ [a:197331, z:[172361], m:2^{31}-1$ \\

\hline
EICG &$ [a:197331, z:[172361], m:2^{48}-59$\\
\hline
\end{tabular}
\end{table}

\begin{table}[!t]
\footnotesize
\caption{Performance of PRNG in CVRP}
\begin{center}
\begin{tabular}{|l|l|l|l|l|l|}
\hline
VRP/PRNG &LehmerCG &LCG &MRG &ICG &EICG\\
\hline
E-n13-k4  &$[257]$ &$[ 257]$ &$[257]$ &$[257]$ &$[ 257]$\\
\hline
E-n22-k4  &$[375]$ &$[ 376]$ &$[390]$ &$[ 384]$ &$[ 381]$\\
\hline
E-n23-k3  &$[610]$ &$[ 583]$ &$[599]$ &$[ 583]$ &$[ 615]$\\
\hline
E-n30-k3  &$[511]$ &$[512]$ &$[526]$ &$[ 513]$ &$[ 508]$\\
\hline
E-n31-k7  &$[471]$ &$[ 470]$ &$[447]$ &$[ 439]$ &$[ 508]$\\
\hline
E-n33-k4  &$[855]$ &$[ 860]$ &$[852]$ &$[ 846]$ &$[ 861]$\\
\hline
E-n51-k5  &$[545]$ &$[ 558]$ &$[587]$ &$[ 574]$ &$[ 566]$\\
\hline
E-n76-k7  &$[713]$ &$[733 ]$ &$[750]$ &$[ 744]$ &$[ 749]$\\
\hline
E-n76-k8  &$[767]$ &$[767 ]$ &$[805]$ &$[ 790]$ &$[809 ]$\\
\hline
E-n76-k10  &$[868]$ &$[889 ]$ &$[850]$ &$[ 879]$ &$[ 889]$\\
\hline

\iffalse
E-n76-k14 &LehmerCG &$[]$\\
\hline
&LCG &$[ ]$\\
\hline
&MRG &$[]$\\
\hline
&ICG  &$[ ]$\\
\hline
%&EICG &$[ 8898.100]$\\
%\hline

\fi

\end{tabular}
\end{center}
\end{table}

In Table VII, we present the performance comparison of Lehmer Congruential Generator, Linear Congruential Generator, Multiple Recursive Generator, Inversive Congruential Generator, and Explicit Inversive Congruential Generator, whose parameters are given in Table VI. 
In the light of the table, we see that Lehmer Congruential Generator and Linear Congruential Generator outperform the others. This could be due to their simplicity and good lattice structure as we reduce each random number modulo a number less than 100 to assign it to a job. On the other hand, we see that Explicit Inversive Congruential Generator has dramatically bad results in almost all cases, and so it shouldn't be used for such applications.

It is widely criticized in the literature that the source code of the implemented algorithms needs to be checked by other researchers \cite{SAC2019}. Besides, sharing the source code in a public domain opens further improvement and this flourishes the research. In this view,  the source code is shared on the web page https://github.com/iedmrc/binary-cws-mcs, which has been implemented according to algorithms in this paper.

\section*{Acknowledgement}
Second and third authors are supported by the Scientific and Technological Research	Council of Turkey (TÜBİTAK) under Project No: \mbox{116R026}.

\end{document}